\documentclass[10pt]{article}
\usepackage{graphicx}
\usepackage{amssymb}
\usepackage{rotating}
\usepackage{pdflscape}
\usepackage[english]{babel}
\usepackage{color}
\usepackage{bigintcalc}
\def\beq{\begin{eqnarray}}
\def\eeq{\end{eqnarray}}
\usepackage{amsfonts}

\usepackage{seqsplit}

\newtheorem{rule-of-thumb}[theorem]{Definition} 

\usepackage{siunitx}

\begin{document}

\title{High order analysis of the limit cycle of the van der Pol oscillator}
\author{Paolo Amore \\
\small Facultad de Ciencias, CUICBAS, Universidad de Colima,\\
\small Bernal D\'{i}az del Castillo 340, Colima, Colima, Mexico \\
\small paolo.amore@gmail.com \\
John P. Boyd \\
\small Department of Atmospheric, Oceanic \& Space Science, \\
\small University of Michigan,  2455 Hayward Avenue, Ann Arbor, MI 48109, USA \\
\small jpboyd@umich.edu \\
Francisco M. Fern\'andez \\
\small INIFTA (CONICET), Divisi\'on Qu\'imica Te\'orica, Blvd. 113 y 64 (S/N), \\
\small Sucursal 4, Casilla de Correo 16, 1900 La Plata, Argentina \\
\small framfer@gmail.com}

\maketitle

\begin{abstract}
We have applied the Lindstedt-Poincar\'e method to study the limit cycle of the van der Pol oscillator,
obtaining the numerical coefficients of the series for the period and for the amplitude to order $859$.
Hermite-Pad\'e approximants have been used to extract the location of the branch cut of the series with unprecendented
accuracy ($100$ digits). Both series have  then been resummed using an approach based on Pad\'e approximants, where the
exact asymptotic behaviors of the period and the amplitude are taken into account. Our results improve drastically all
previous results obtained on this subject.
\end{abstract}

\section{Introduction}
\label{sec_intro}

The van der Pol equation is possibly the most known example of ordinary differential equation with a limit
cycle; it reads
\begin{eqnarray}
\ddot{x}(t) + x(t)  = \mu \dot{x}(t) (1-x(t)^2)
\label{eq_vdp}
\end{eqnarray}
where $\mu>0$ is a parameter which controls the strength of the
nonlinear contributions. Even though for $\mu = 0$ this equation
reduces to a simple harmonic oscillator, that can sustain
oscillations of arbitrary amplitude, for any $\mu >0$ the van der
Pol oscillator has a single oscillatory mode, with a specific
amplitude and period, which depends on $\mu$. This mode is usually
referred to as a ``limit cycle". When the initial conditions are
far from the limit cycle, the solution decays exponentially fast
towards the limit cycle; this behavior results into a trajectory
in phase space which spirals out (in) if the initial point falls
inside (outside) the limit cycle, rapidly reaching  the closed
trajectory.

The reader may find a historical account of the origins of the van
der Pol equation in a recent paper by Ginoux and Letellier,
ref.~\cite{Ginoux12}; for a general perspective on the subject of
self-sustained oscillations, it is useful to refer to the recent
review article by Jenkins~\cite{Jenkins13}. Self-oscillations, as
those described by Eq.~(\ref{eq_vdp}),  cover an important role in
different areas of Physics and Mathematics, and even of Biology:
van der Pol and van der Mark, for example, described the heartbeat
in terms of a relaxation oscillation and provided an electrical
model of the heart~\cite{VanderPol28}; FitzHugh,
\cite{FitzHugh61}, has generalized van der Pol equation to
describe the excitation of neurons and muscle fibers.

The task of finding approximate solutions to Eq.~(\ref{eq_vdp})
becomes challenging, particularly when $\mu$ is large, because the
system displays two alternating slow and fast regions, which
reflect in an highly deformed cycle in phase space. We may
classify the different strategies that have been used in the
literature to solve the van der Pol equation  into three main
categories: numerical studies, carried out up to some large but
finite value of $\mu$, asymptotic calculations, valid for $\mu
\rightarrow \infty$, and  perturbative calculations, obtained
expanding about $\mu=0$.

As nowadays the numerical solution of Eq.~(\ref{eq_vdp}) is easily
accessible due to the wide disponibility of hardware and software,
we mention the works in the first category,
refs.~\cite{Krogdahl60,Urabe60,Clenshaw66,Greenspan72}, mainly for
historical reasons.

The first example of work belonging to the second class is given by refs.~\cite{Haag43}, where a formula for the asymptotic period
was given (as reported by Ginoux and Letellier \cite{Ginoux12})
\begin{eqnarray}
T^{(asym)}(\mu) = (3-\log 4) \mu + \frac{12.89}{\mu^{1/3}} + \frac{2}{\mu} \left(-3.31 + \frac{19}{9} \log \mu \right) - \frac{4}{\mu^{5/3}} + \dots
\end{eqnarray}

More precise formulas have then been derived in \cite{Dor52, Urabe60, Ponzo65, Zonneveld66}; the leading terms in the expression
for the period read
\begin{eqnarray}
T^{(asym)}(\mu) = (3-\log 4) \mu - \frac{3\alpha}{\mu^{1/3}} - \frac{2}{3 \mu} \log\mu - \frac{1.3246}{\mu^2}  + \dots
\label{eq_asym_T}
\end{eqnarray}
where $\alpha \approx -2.33810741$ is the highest zero of the Airy function.

An analogous formula for the amplitude of the oscillations has also been derived in \cite{Dor52, Urabe60, Ponzo65} and it reads
\begin{eqnarray}
A^{(asym)}(\mu) = 2   - \frac{\alpha}{3 \mu^{4/3}} - \frac{16}{27} \frac{\log \mu}{\mu^2}
+   \frac{\left(3\beta -1+\log 4 -8 \log 3 \right)}{9 \mu^2} + O(\mu^{-8/3})
\label{eq_asym_A}
\end{eqnarray}
where $\beta \approx 0.1723$.

Finally the last class has to do with the perturbative calculation
of the solutions to Eq.~(\ref{eq_vdp}): in this class fall the
papers by Deprit and Rom ~\cite{Deprit67,Deprit68} and Deprit and
Schmidt \cite{Deprit79}, by Andersen and Geer \cite{Andersen82}
and Dadfar, Andersen and Geer \cite{Dadfar84}, by Buonomo
\cite{Buonomo98b} and, more recently, by Suetin \cite{Suetin12}.
The Lindstedt-Poincar\'e method allows to obtain expressions for
the period and the amplitude as power series in $\mu$; the
presence of a finite radius of convergence of these series,
however, limits the direct application of the perturbative
formulas to the domain of convergence. Different resummation
procedures have been discussed in
refs.~\cite{Andersen82,Dadfar84,Buonomo98b,Suetin12}, finding a
radius of convergence $\mathcal{R} \approx 3.42$. In this class we
also mention ref.\cite{Amore04}, based on an alternative
implementation of the Lindstedt-Poincar\'e method, with improved
convergence properties.

In the present work we have set the following goals:
\begin{itemize}
\item To obtain the largest number of coefficients of the perturbative series for the amplitude and the period of the van der Pol oscillator;
\item To obtain the most precise determination of the branch cut of these series;
\item To perform a non-perturbative resummation of the power series, taking advantage of the goals previously achieved;
\end{itemize}

The paper is organized as follows: in section \ref{sec_LP} we describe the Lindstedt-Poincar\'e method and compare our result with previously obtained
results; in section \ref{sec_HP} we apply the Hermite-Pad\'e approximants to obtain a determination of an accurate location of the branch cut;
the resummation of the series is carried out in section \ref{sec_resum}; finally the conclusions are stated in section \ref{sec_concl}.

\section{Lindstedt-Poincar\'e method}
\label{sec_LP}

We consider the van der Pol equation (\ref{eq_vdp}). As we have mentioned in the Introduction,
this equation has a limit cycle, meaning that for $t \rightarrow \infty$ the solutions tend to
a periodic solution with fixed amplitude, regardless of the initial conditions.
In a straightforward application of perturbation theory, where the solution can be expressed as a
power series in $\mu$, secular terms appear, which completely spoil the expansion.

In the Lindstedt-Poincar\'e method one introduces a ``universal"
time $\tau = \Omega(\mu) t$ (``strained coordinate"), where
$\Omega(\mu)$ is the exact frequency of the oscillations; the
coefficients of the power series in $\mu$ for $\Omega$ are
determined requiring that the ``secular terms" arising at each
perturbative order are avoided. In the particular case of the van
der Pol oscillator, and more in general of problems with a limit
cycle, a similar power series in $\mu$ must be assumed for the
amplitude.

Let us briefly discuss in detail the implementation of the method for the van der Pol equation.
After switching to the new time variable, we define $y(\tau) \equiv x(t)$
and the original differential equation takes the form
\begin{eqnarray}
\Omega^2 \frac{d^2y}{d\tau^2} + y(\tau) - \mu \Omega \frac{dy}{d\tau} \left( 1 - y^2(\tau) \right) = 0
\label{eq_vdp_omega}
\end{eqnarray}
with the initial conditions
\begin{equation}
y(0) = \mathcal{A}(\mu) \ \ \ ; \ \ \ \frac{dy}{d\tau}(0) = 0 \nonumber
\end{equation}

In our perturbative scheme we will assume that
\begin{eqnarray}
\Omega(\mu) = 1 + \sum_{n=1}^\infty \omega_n \mu^n
\label{eqomega}
\end{eqnarray}
and
\begin{eqnarray}
\mathcal{A}(\mu) = \sum_{n=0}^\infty a_n \mu^n
\label{eqA}
\end{eqnarray}
where $\omega_n$ and $a_n$ are constant coefficients and $\omega_0=1$ is the frequency of the linear system.

Similarly we look for a solution
\begin{eqnarray}
y(\tau) = \sum_{n=1}^\infty y_n(\tau) \mu^n \nonumber
\end{eqnarray}
obeying the initial conditions
\begin{eqnarray}
y_n(0) = a_n \ \ \ ; \ \ \ \frac{dy_n}{d\tau}(0)= 0 \nonumber
\end{eqnarray}

Upon substitution of these expressions in the differential equation,
one has an infinite set of coupled {\sl linear} non--homogeneous differential
equations, corresponding to the different orders in $\mu$, which
can be solved sequentially starting from the lowest order.

For instance to order $\mu^0$ one has
\begin{eqnarray}
\frac{d^2 y_0}{d\tau^2} + y_0(\tau) = 0 \nonumber
\end{eqnarray}
with solution
\begin{eqnarray}
y_0(\tau) = a_0 \cos\tau \nonumber
\end{eqnarray}

To order $\mu$, using the solution of order $0$, one has the differential equation
\begin{eqnarray}
\frac{d^2 y_1}{d\tau^2} + y_1(\tau) 
&=& \frac{a_0 (a_0^2-4)}{4} \sin\tau +  2 a_0 \omega_1 \cos\tau + \frac{a_0^3}{4} \sin 3\tau  \nonumber
\end{eqnarray}
where in the last line the resonant terms containing $\sin  \tau$
and $\cos  \tau$ are responsible of the appearance of ``secular
terms" and therefore they must be avoided. This condition requires
\begin{eqnarray}
a_0 = 2 \ \ \ ; \ \ \ \omega_1 = 0  \nonumber
\end{eqnarray}

The solution to the differential equation then takes the form
\begin{eqnarray}
y_1(\tau) = a_1 \cos \tau + \frac{3}{4} \sin \tau - \frac{1}{4} \sin 3\tau \nonumber
\end{eqnarray}
Note that the first order solution contains the higher frequency
$3$.

The Lindstedt-Poincar\'e method can be efficiently implemented on a computer, to very large orders,
to calculate the power series expansions for $\Omega$, $\mathcal{A}$ and $y$, either exactly or
numerically. We have written two different programs, one for the symbolic calculation of the power
series, and one for the numerical calculation, that has been performed working  with $1000$ digits.

In Table \ref{Table_orders_LP} we compare the orders of the LP method obtained in previous works with the
orders obtained in the present work. The disposal of a large number of coefficients of the perturbative
series, combined with the high accuracy of the calculation, will be essential  for the numerical
analysis performed in the next sections.

\begin{figure}[t]
\begin{center}
\includegraphics[width=8cm]{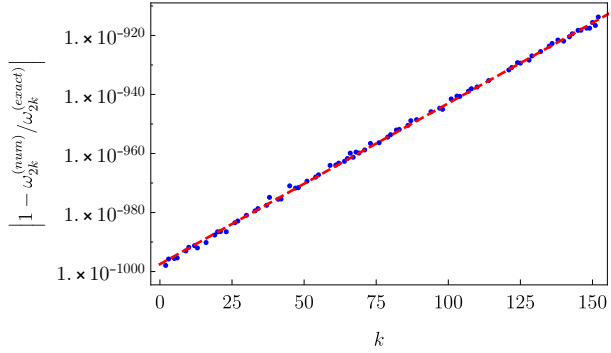}
\caption{Relative error over the numerical coefficients of the order $2k$, with respect to the exact ones. The dashed line is the fit $10^{997.5-0.545 \ x}$.}
\label{fig_1}
\end{center}
\end{figure}

\begin{table}
\begin{center}
\begin{tabular}{|c|c|c|}
\hline\noalign{\smallskip}
& {\rm order} & {\rm order} \\
& [{\rm symbolic}] & [{\rm numerical}] \\
\noalign{\smallskip}\hline\noalign{\smallskip}
{\rm Deprit and Rom,}  \ {\rm Ref.} \cite{Deprit67}         &   -  & $30$  ($12$ digits)\\
{\rm Deprit and Schmidt,}  \ {\rm Ref.} \cite{Deprit79}     &   8  & - \\
{\rm Andersen and Geer,}  \ {\rm Ref.} \cite{Andersen82}    &  24  & $163$ ($37$ digits)\\
{\rm Buonomo,}  \ {\rm Ref.} \cite{Buonomo98b}              & 120  & $402$ ($40$ digits) \\
{\rm Suetin,}   \ {\rm Ref.} \cite{Suetin12}                &   -  & $442$ ($40$ digits) \\
{\rm This work}                                             & 308  & $859$ ($1000$ digits) \\
\noalign{\smallskip}\hline
\end{tabular}
\caption{Comparison of the orders of the perturbative expansion for the van der Pol oscillator, obtained with the
Lindstedt-Poincar\'e method on a computer, in different works. }
\label{Table_orders_LP}
\end{center}
\end{table}

To assess the quality of the numerical coefficients, in Fig.~\ref{fig_1} we have plotted the quantity $-\log_{10} \left| 1 - \frac{\omega_{2j}^{(num)}}{\omega_{2j}^{(exact)}}\right|$, as a function of the order of the coefficients themselves.
The dashed line in the plot is the linear fit
\begin{equation}
f(x) = 997.5-0.273\ x
\end{equation}
which vanishes at $x \approx 3660$ (thus $n_{max} \approx 3660$ is approximately the highest order of the Lindstedt-Poincar\'e
expansion that can be obtained working with $1000$ digits).



\section{Hermite-Pad\'e approximants}
\label{sec_HP}

Applying the Lindstedt-Poincar\'e method as explained in the
previous section one is able to obtain  a large number of terms of
the perturbative series for the period and for the amplitude of
the van der Pol oscillator. Because these series are powers  in
the \emph{square} of $\mu$, it is convenient to define the new
expansion parameter $\nu \equiv \mu^2$. The period and amplitude
series then take the form
\begin{eqnarray}
U_N = \sum_{n=0}^{N-1} c_n \nu^n
\end{eqnarray}
where the coefficients $c_n$ are either calculated {\sl exactly} or numerically, with high precision.

These series converge for $|\nu| < 3.42$ [$|\mu| < 1.85$] and diverge for larger $|\nu|$ where the radius of
convergence is computed with high accuracy below.

Since arbitrarily large values of $\nu$ are of interest and calculations for $\nu$ as large as 10,000 have been previously published,
one is then faced with a task of extracting useful  information from a series far beyond its radius of convergence.

Suetin applied Pad\'e approximation
\cite{Suetin12} to analyze the singularities of the amplitude and frequency functions, building on work in
\cite{Andersen82,Dadfar84}. He conjectured (``it is possible that the following assertion is valid"):

\begin{enumerate}
\item The singularities of the amplitude and period functions are of the form
\begin{eqnarray}
g(\nu) \, \left\{ \exp(i \alpha) \sqrt{ \nu - \mathcal{R} \exp(i \phi) } + \exp(-i \alpha)
\sqrt{ \nu - \mathcal{R} } \right\}
\end{eqnarray}
where $g(\nu)$ is holomorphic in a disk whose diameter is larger than
$\mathcal{R}$.
\item The radius of convergence $\mathcal{R}$ is exactly 3.42.
\item The phase $\phi$ is exactly $89 \pi/156$.
\item $\alpha$ is $\pi/4$ for the frequency and $3 \pi/10$ for the amplitude.
\end{enumerate}

However, Pad\'e approximants  are always rational functions, and therefore are a clumsy tool
for investigating square-roots. We therefore use a generalization of Pad\'e approximations,
 {\it Hermite-Pad\'e}  approximants, which include square root singularities.  Hermite-Pad\'e  approximations
have been successfully used in rootfinding \cite{Shafer74} and pgs. 46-47 of \cite{BoydBook3},
liquid drops \cite{ErRianiElJarroudiSeroGuillaume14} and water waves \cite{Common82},
quantum mechanics~\cite{VainbergMurPopovSergeev86,GermannKais93,BoydOP118,
LopezCabreraGoodsonHerschbachMorgan92,CabayJonesLabahn97,Sergeev95,SuvernevGoodson97,
SergeevGoodson98,MayerTong85,MayerNuttallTong84,Gammel73,ErRianiElJarroudiSeroGuillaume14}, series analysis of multivalued functions \cite{DrazinTourigny96} and other applications
\cite{Shafer74,VainbergMurPopovSergeev86,GermannKais93,BoydOP118,DrazinTourigny96,
LopezCabreraGoodsonHerschbachMorgan92,CabayJonesLabahn97,Sergeev95,SuvernevGoodson97,SergeevGoodson98,MayerTong85,MayerNuttallTong84,Gammel73,ErRianiElJarroudiSeroGuillaume14}.
The Cauchy rootfinding method employs a low order quadratic approximation~\cite{BoydBook3}.
A Hermite- Pad\'{e} approximation is defined as the solution to a polynomial equation of arbitrary user-chosen degree. Here,
we use only approximations that solve a quadratic equation.The name ``Shafer approximation" is also used as a
synonym for quadratic-solving Hermite-Pad\'{e} approximations \cite{Shafer74}. Here,we use only Shafer approximations.


Many illuminating and useful theorems have been proven about convergence properties of the Hermite-Pad\'e approximants.
However, numerical ill-conditioning is not unusual, so care must be taken and sometimes
multiple precision arithmetic is helpful in achieving high accuracy.


To approximate a function $f(x)$ which is known only through the first $N$ term of its power series, use the function $f_[K/L/M]$,
the quadratic Shafer approximation, which is defined to be the solution of the quadratic equation
\begin{equation}
P (f[K/L/M])^{2} + Q \, f[K/L/M] \, + \, R = 0
\end{equation}
where the polynomials $P$, $Q$ and $R$ are of degrees $K$, $L$ and $M$, respectively.
These polynomials are chosen so that the power series expansion of $f[K/L/M]$ agrees with that of $f$ through
the first $N=K+L+M+1$ terms.

The constant terms in $P$ and $Q$ can be set arbitrarily without
loss of  generality since these choices do not alter the root of
the equation, so the total number of degrees of freedom is as
indicated. As for ordinary Pad\'e approximants, empirically, the
most accurate approximations are usually obtained by choosing the
polynomials to be of equal degree, the so-called ``diagonal"
approximants.

Again as for ordinary Pad\'e approximants, the coefficients of the
polynomials can be computed by solving a matrix equation. Cabay, Jones and Labahn \cite{CabayJonesLabahn97}
describe a faster but more ill-conditioned recursive algorithm to compute
 Hermite-Pad\'e approximants. Mayer, Nuttall and Tong present and apply a simple recursion
restricted to the quadratic Shafer approximation \cite{MayerNuttallTong84,MayerTong85}. Loi and McInnes
also offer an algorithm to compute  the quadratic Hermite-Pad\'e approximation \cite{LoiMcInnes84}

Like  Pad\'e approximants, Hermite-Pad\'e approximants in a parameter $\nu$ are \emph{highly nonuniform} in $\nu$. The approximations are extremely accurate for small
$\nu$, but error grows exponentially with increasing $|\nu|$. Nevertheless, Hermite-Pad\'e approximations are usually convergent far beyond the radius of convergence of the power series from which they were constructed. Indeed, they may converge exponentially fast for all
finite $\nu$ even when the power series is factorially divergent for all non-zero $\nu$.

We applied classical Pad\'{e} methods and stopped where these met
our needs. When Pad\'{e} methods disappointed, we shifted to
Hermite-Pad\'{e}.

Parenthetically, note that these power series based algorithms can be extended to
Chebyshev interpolation as in the Hermite-Pad\'e-Chebyshev work of Boyd \cite{BoydOP128}.
 There is an extensive literature on Chebyshev polynomial generalizations of ordinary
Pad\'e approximations \cite{MasonCrampton05,MasonCrampton05b}. The Chebyshev-based
generalizations are much more uniform in accuracy on the interval spanned by the
interpolation points than power series-based accelerations, but cannot be applied
directly to perturbation series. Rather, the Chebyshev-Pad\'e algorithms are useful for
converting numerical solutions at a set of discrete values of a parameter $\nu$ into
 spectrally-accurate explicit approximations that are continuous in $\nu$.

\begin{figure}[t]
\begin{center}
\includegraphics[width=8cm]{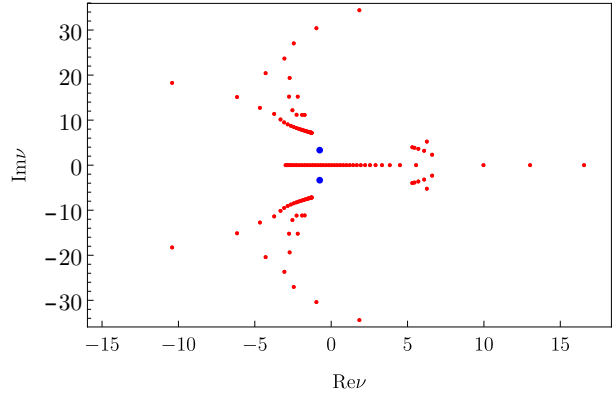}
\caption{Roots of $Q^2-4 P R$ for the case $K=L=M=142$; the blue ones correspond to the roots in the table.}
\label{fig_HP_2}
\end{center}
\end{figure}

In the complex $\nu$ plane there are two complex conjugate roots,
\begin{eqnarray}
\nu_{\pm} = \mathcal{R} \ e^{\pm i \phi}
\end{eqnarray}
at which $\Delta$ vanishes.  The blue points in the figure \ref{fig_HP_2} correspond to these complex roots
for the case $K=L=M=142$.

\begin{figure}[t]
\begin{center}
\includegraphics[width=6cm]{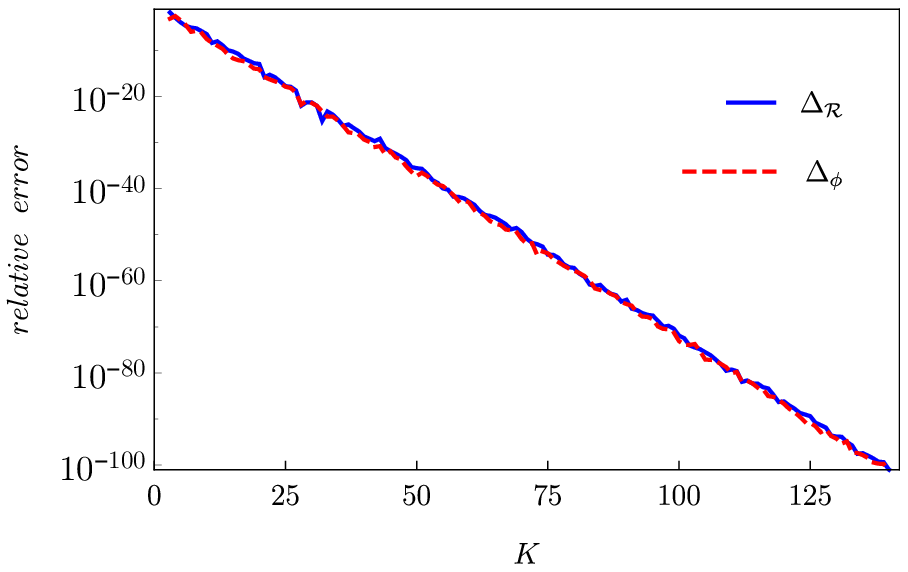} \hspace{.5cm}
\includegraphics[width=6cm]{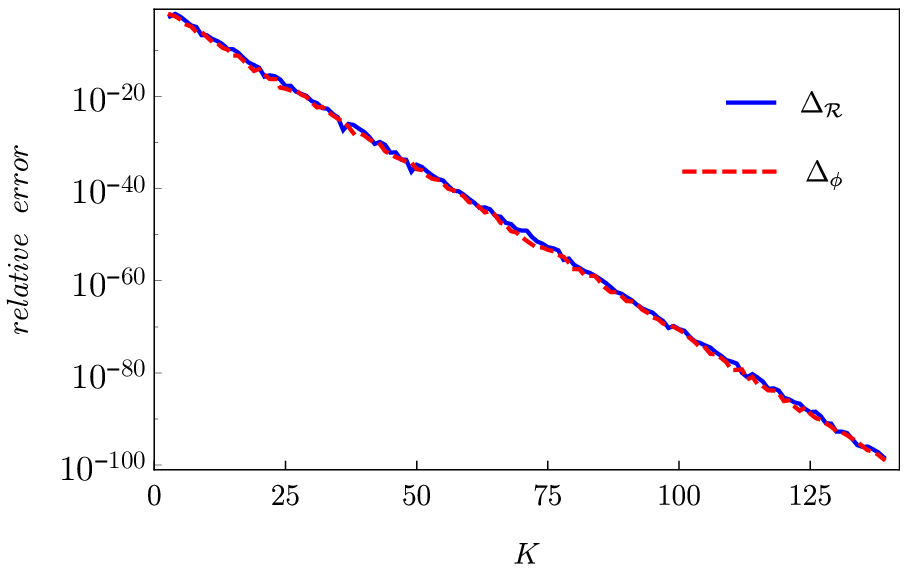}
\caption{Relative errors over the radius of convergence $\mathcal{R}$ and the phase $\phi$ estimated with the Hermite-Pad\'e method for the series
of the period (left plot) and for the series of the amplitude (right plot). Here $\Delta_{\mathcal{R}} \equiv | \mathcal{R}^{(K)}-\mathcal{R}^{(142)}|$
and $\Delta_{\phi} \equiv | \mathcal{\phi}^{(K)}-\mathcal{\phi}^{(142)}|$. }
\label{fig_HP_1}
\end{center}
\end{figure}

\begin{figure}[t]
\begin{center}
\includegraphics[width=8cm]{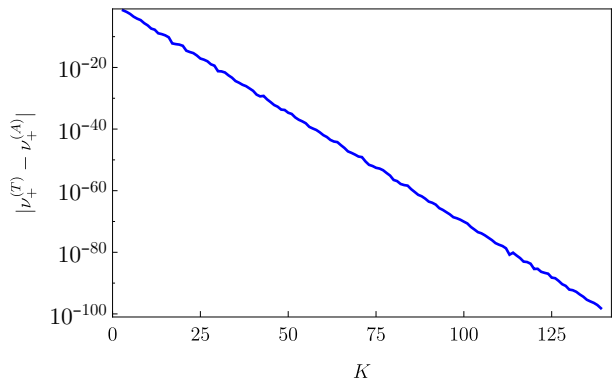}
\caption{$| \nu_+^{(T)}-\nu_+^{(A)}|$ as a function of the order $K$.}
\label{fig_HP_3}
\end{center}
\end{figure}

We have calculated the diagonal Hermite-Pad\'e approximants up to a maximal order $K=142$, which corresponds to using $854$
terms of the numerical series (where the coefficients of the series have been calculated with $1000$ digits of accuracy).
We have performed our analysis for both the series for the period and for the series of the amplitude: the plots in
Figs.\ref{fig_HP_1}  illustrate  the error over $\mathcal{R}$ and $\phi$ for both series, as a function  of the order
of the approximant. In Fig.~\ref{fig_HP_3}  we plot the difference $|\mathcal{R}_T-\mathcal{R}_A|$, for the radius of convergence
of the two series. The two results agree to at least $100$ digits!

Our most accurate results are (where all digits are expected to be be accurate)
\begin{eqnarray}
\mathcal{R} &=& \num{3.420187909357135029477567375233640987583555906152748024797507314} \nonumber \\
&& \num{779205532661443624837516502144919691} \nonumber  \\
       \phi &=& \num{1.792288671545795214263603138015353292564791705087259065160879953} \nonumber \\
&& \num{197027910943997296402573377464202939} \nonumber
\end{eqnarray}
 These results falsify both Suetin's conecture that $\mathcal{R}$ is exactly 3.42 and
his suggestion that $\phi=89 \pi/156$, which is in error by 0.0000302.

Suetin mentions on pg. S25 that he applied diagonal Hermite-Pad\'e approximants to
convince himself that the character of the singularities of the frequency and amplitude
are more complicated than quadratic branching [square root singularities],
but we disagree.
A straightforward application of the generating function method developed by Fern\'andez et al.~\cite{Fernandez87}
that is a simplification of the approach proposed earlier by Hunter and Guerrieri ~\cite{Hunter80} confirms that the singularity
is of square-root type.

\begin{figure}[t]
\begin{center}
\includegraphics[width=8cm]{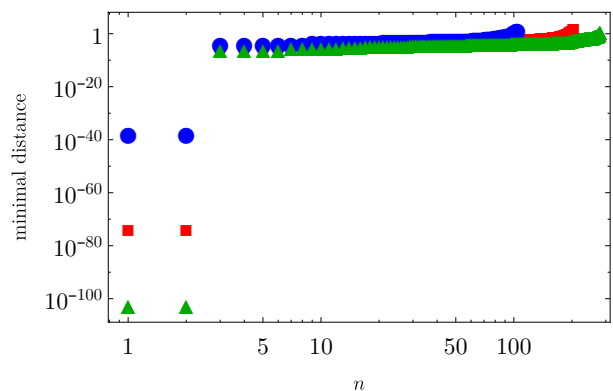}
\caption{Minimal distance $d_n(K)$ of each of the roots of the Hermite-Pad\'e approximant of order $K$ ($n$ is the index counting these roots)  from any of the roots of order $K+1$, for $K=50$ (circles), $K=100$ (squares) and $K=138$ (triangles).}
\label{fig_HP_min_dist}
\end{center}
\end{figure}

The radical has many zeros besides the pair shown in blue in
Fig.~\ref{fig_HP_2}. The many spurious zeros can be identified by
comparing calculations with different $K$. The relevant zeros of
$R(\nu)$ are convergent whereas the spurious zeros hop around as
the order $K$ of the   Pad\'e approximation is varied. This
behavior is illustrated in Fig.~\ref{fig_HP_min_dist}, where we
have plotted the minimal distance $d_n(K)$ of each of the roots of
the Hermite-Pad\'e approximant of order $K$ from any of the roots
of order $K+1$, for $K=50$ (solid line), $K=100$ (dashed line) and
$K=138$ (dotted line). It is evident that, while two of the roots
converge to a given value (these are precisely the roots that we
have mentioned earlier), the remaining roots to not converge.

\section{Resummation}
\label{sec_resum}

The precise determination of the branch points of the series for the period and for the amplitude can be exploited to obtain an effective resummation
of the series itself.

To this purpose it is convenient to introduce  the variable $\xi(\nu) \equiv (\nu^2 -2\nu \mathcal{R} \cos\phi + \mathcal{R}^2)^{1/4}$, which obeys the limit
\begin{eqnarray}
\lim_{\nu \rightarrow \infty} \frac{\xi(\nu)}{\sqrt{\nu}} = 1 \nonumber
\end{eqnarray}

On the other hand, the leading asymptotic behavior of the period,
for $\mu \rightarrow \infty$, has been obtained by Dorodnitsyn
\cite{Dor52} (see Eq.~(\ref{eq_asym_T}).

A simple strategy consists of ``depurating" the period of the
leading asymptotic behavior given above by introducing
\begin{eqnarray}
\tilde{T}(\nu) &\equiv&  \frac{2\pi}{\Omega} -(3 - \log 4) \xi + \frac{3\alpha}{\xi^{1/3}} + \frac{2 \log \xi}{3\xi} \nonumber \\
&+& \left( 1.3246 + \frac{1}{2} \mathcal{R}\cos\phi (-3 + \log 4) \right) \frac{1}{\xi}
\end{eqnarray}
and then obtain the ``complete" period as
\begin{eqnarray}
T^{(resummed)} &=& \left[ \tilde{T}(\nu)\right]_{M,N} +(3 - \log 4) \xi - \frac{3\alpha}{\xi^{1/3}} - \frac{2 \log \xi}{3\xi} \nonumber \\
&-& \left(1.3246 + \frac{1}{2}  \mathcal{R} \cos\phi (-3 + \log 4) \right) \frac{1}{\xi}
\end{eqnarray}
where $\left[ \tilde{T}\right]_{M,N}$ is the $[M,N]$ Pad\'e approximant to $\tilde{T}$. In the present case, since the asymptotic behavior up to order $1/\mu$ has been taken care of,
it is convenient to use $N=M+1$.

We can go a step further and try to guess the next contribution in
the asymptotic formula; it makes sense to assume that this term
behaves as $\gamma/\mu^{4/3}$, where $\gamma$ is a parameter, and
therefore we write the ``depurated" period as
\begin{eqnarray}
\tilde{T}(\nu) &\equiv&  \frac{2\pi}{\Omega} -(3 - \log 4) \xi + \frac{3\alpha}{\xi^{1/3}} + \frac{2 \log \xi}{3\xi} \nonumber \\
&+& \left( 1.3246 + \frac{1}{2} \mathcal{R}\cos\phi (-3 + \log 4) \right) \frac{1}{\xi} - \frac{\gamma}{\xi^{4/3}}
\end{eqnarray}

As before, we obtain the ``complete" period as
\begin{eqnarray}
T^{(resum)}(\gamma) &=& \left[ \tilde{T}(\nu)\right]_{M,N} +(3 - \log 4) \xi - \frac{3\alpha}{\xi^{1/3}} - \frac{2 \log \xi}{3\xi} \nonumber \\
&-& \left(1.3246 + \frac{1}{2} \mathcal{R}\cos\phi (-3 + \log 4) \right) \frac{1}{\xi} +  \frac{\gamma}{\xi^{4/3}}
\end{eqnarray}

\begin{figure}[t]
\begin{center}
\includegraphics[width=8cm]{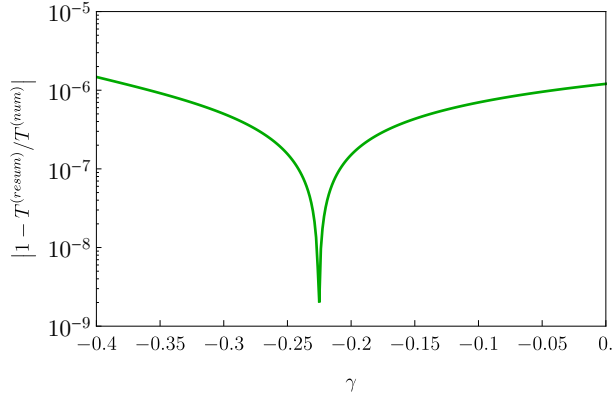}
\caption{Relative error $\left|1 - T^{(resum)}/T^{(num)}\right|$ at $\mu = 100$, as a function of $\gamma$.}
\label{fig_HP_5}
\end{center}
\end{figure}

In Fig.~\ref{fig_HP_5} we plot the relative error $\left|1 - T^{(resummed)}/T^{(num)}\right|$ at $\mu = 100$, as a function of $\gamma$; the
Pad\'e approximant in the resummed expression is calculated using $M=N-1=120$. For $\gamma = 0$ one recovers the results obtained earlier with the
standard asymptotic formula (\ref{eq_asym_T}).  The error is minimal for $\gamma \approx - 0.225$ (roughly two orders of magnitude better!).

\begin{figure}[t]
\begin{center}
\includegraphics[width=8cm]{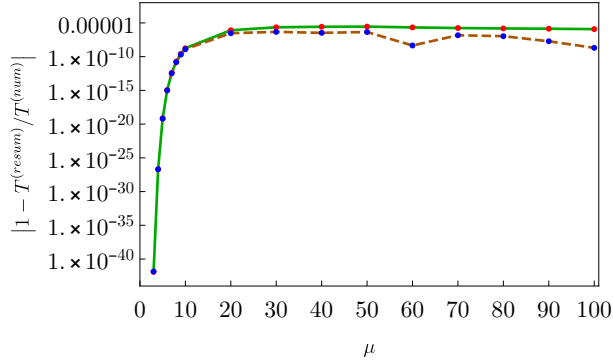}
\caption{
Relative error over the period as a function of $\mu$ using $T^{(resum)}(0)$ (solid line) and $T^{(resum)}(-0.225)$ (dashed line). The resummed expression has been calculated using a Pad\'e (not Hermite-Pad\'e) with $M=N-1=100$.}
\label{fig_HP_4}
\end{center}
\end{figure}

The error $\left|1 - T^{(resummed)}/T^{(num)}\right|$ is plotted in Fig.~\ref{fig_HP_4},  using $T^{(resum)}(0)$ (solid line) and $T^{(resum)}(-0.225)$ (dashed line). The resummed expression has been calculated using a Pad\'e with $M=N-1=120$.

Let us now turn our attention to the amplitude.  Eq.~(\ref{eq_asym_A}) is the asymptotic expression for the amplitude for $\mu \rightarrow \infty$
derived by Dorodnitsyn and also confirmed later by Urabe \cite{Urabe60} and by Ponzo and Wax \cite{Ponzo65}.

\begin{figure}[t]
\begin{center}
\includegraphics[width=8cm]{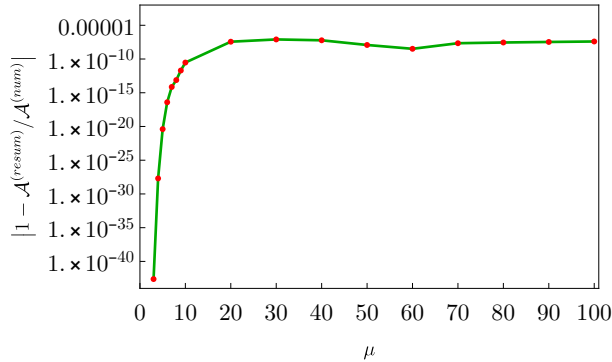}
\caption{Relative errors over the amplitude  as a function of $\mu$ using $\mathcal{A}^{(resum)}$. The resummed expression has been calculated
using a Pad\'e (not Hermite-Pad\'e)  with $M=N-2=100$.}
\label{fig_HP_6}
\end{center}
\end{figure}

As for the case of the period we introduce a quantity where the leading asymptotic behaviors are eliminated
\begin{eqnarray}
\tilde{\mathcal{A}}(\nu) &\equiv&  \mathcal{A} -2 +  \frac{\alpha}{3 \xi^{4/3}} - \frac{16}{27} \frac{\log \xi}{\xi^2}
+ \frac{1}{9 \xi^2}  \left(3\beta -1+\log 4 -8 \log 3 \right)
\end{eqnarray}
and then obtain the ``complete" amplitude as
\begin{eqnarray}
\mathcal{A}^{(resum)} &=& \left[ \tilde{\mathcal{A}}(\nu)\right]_{M,N} +2 -  \frac{\alpha}{3 \xi^{4/3}}
+ \frac{16}{27} \frac{\log \xi}{\xi^2} \nonumber \\
&-& \frac{1}{9 \xi^2}  \left(3\beta -1+\log 4 -8 \log 3 \right)
\end{eqnarray}
where $ \left[ \tilde{\mathcal{A}}(\nu)\right]_{M,N}$ is the Pad\'e approximant of orders $M$ and $N$. In the present case, it is convenient
to choose $M=N-2$, given that the asymptotic formula (\ref{eq_asym_A}) contains a term with the behavior $1/\mu$.
The relative error on the amplitude is shown in the  plot of Fig.~\ref{fig_HP_6}.

\section{Conclusions}
\label{sec_concl}

In this paper we have applied the Lindstedt-Poincar\'e method to the van der Pol equation, to high orders in the expansion parameter.
We have been able to reach order $308$ in a purely symbolic calculation, where the perturbative coefficients are calculated in arbitrary
precision, and order $859$ in a numerical calculation, where the perturbative coefficients are calculated with an accuracy of $1000$
digits. A comparison with similar calculations in the literature shows that our results are the most precise -- most orders, most
digits --  (see Table \ref{Table_orders_LP}).

Using the precise coefficients both for the frequency (period) and for the amplitude, we have been able to estimate with unprecedented
accuracy the radius of convergence $\mathcal{R}$ of the series and the location of the branch cut, using the Hermite-Pad\'e approximant.
In all the previous estimations, the radius of convergence and the phase ($\mathcal{R}$ and  $\phi$ respectively)  had been obtained
with few digits; for example, Andersen and Geer \cite{Andersen82}, report $\mathcal{R} \approx 3.42$ and $\phi \approx 1.7925$, and a similar accuracy
is also achieved by Suetin \cite{Suetin12}. Having at our disposal more (and more precisely calculated) series coefficients and
using a different technique, the Hermite-Pad\'e approximants, we have obtained results for $\mathcal{R}$ and $\phi$ which are expected to have the first
$100$ digits correct. It is remarkable that the numerical value of the branch cut for the series for the period and for the amplitude
converge to the same result, within $100$ digits.

This extremely accurate determination of the brach cut can be used to obtain a better resummation of the perturbative series. Our resummation
involves a mapping of the perturbative parameter $\mu$ to $\xi(\nu) = (\nu^2 -2\nu \mathcal{R} \cos\phi + \mathcal{R}^2)^{1/4}$ (similar to the one used
by Andersen and Geer), followed by a straightforward application of Pad\'e approximants, taking into account the asymptotic behavior
of the period and of the amplitude. The resummed expressions that we have obtained are much more accurate than similar expressions
previously found (for a comparison see for instance \cite{Andersen82,Buonomo98b}), over all range of $\mu$ considered.

\section*{Acknowledgements}
The research of P.A. was supported by the Sistema Nacional de Investigadores (M\'exico).
J.P. Boyd was supported by the U. S. National Science Foundation through grant DMS-1521158.

\end{document}